\documentclass[12pt, a4paper]{amsart} 

\usepackage{amsfonts}
\usepackage{epsfig}

\newcommand{\dd}{\;\mathrm{d}}

\newenvironment{ack}{\noindent{\small{A}\tiny{CKNOWLEDGMENTS}.}}{} 
\newenvironment{keyword}{\noindent{\small{K}\tiny{EYWORDS}.}}{} 
\newenvironment{pf*}[1]{\noindent {\it #1}}{$\Box$} 
\newenvironment{pf}{\noindent {\it Proof.}}{$\Box$} 
\newcommand{\e}{\mathrm{e}}
\newcommand{\Rset}{\mathbb{R}} 
\newcommand{\Cset}{\mathbb{C}} 
\newcommand{\Qset}{\mathbb{Q}} 
\newcommand{\Nset}{\mathbb{N}} 
\newcommand{\Zset}{\mathbb{Z}} 
\newtheorem{thm}{Theorem} 
\newtheorem{defn}[thm]{Definition} 
\newtheorem{lem}[thm]{Lemma} 

\newtheorem{theorem}{Theorem}[section]
\newtheorem{lemma}[theorem]{Lemma}
\newcommand{\A}{\mathcal A}
\newcommand{\RA}{{\mathrm{Res}}_{\A_0,\ldots,\A_n}}
\newcommand{\Res}{\mathrm{Res}}
\newcommand{\TT}{\mathbb{T}}
\newcommand{\PP}{\mathbb{P}}
\newcommand{\Li}{\mathrm{Li}} 
\newcommand{\re}{\mathop{\mathrm{Re}}} 
\newcommand{\im}{\mathop{\mathrm{Im}}} 
\newcommand{\ii}{\mathrm{i}} 
\newcommand{\Lf}{\mathrm{L}} 

\begin{document}


\title{On the Mahler measure of resultants in small dimensions}
\author{Carlos D'Andrea}
\address{Department d'\'Algebra i Geometr\'ia, Universitat de Barcelona, Gran Via de les Corts Catalanes 585, 08007 Barcelona. Spain}
\email{carlos@dandrea.name}
\author{Matilde N. Lal\'{\i}n}
\address{Mathematical Sciences Research Institute,
17 Gauss Way,
Berkeley, CA 94720}
\email{mlalin@math.ubc.edu}

\begin{abstract} We prove that sparse resultants having Mahler measure
equal to zero are those whose Newton polytope has dimension one. We then
compute the Mahler measure of resultants in dimension two, and examples in
dimension three and four. Finally, we show that sparse resultants are
tempered polynomials. This property suggests that their Mahler measure may lead to special
values of $\Lf$-functions and polylogarithms. \end{abstract}
\maketitle 

\begin{keyword}
Sparse resultants, Mahler measure, height, Newton polytope, polylogarithms, tempered polynomials. 
\end{keyword}


\section{Introduction}
Let $\A_0,\ldots,\A_n\subset\Zset^n$ be finite sets of integral
vectors, $\A_i:=\{a_{ij}\}_{j=1,\ldots,k_i}.$ We denote with
$\RA\in\Zset[X_0,\ldots,X_n]$ the associated mixed sparse
resultant, which is an irreducible polynomial in $n+1$ groups
$X_i:=\{x_{ij};\, 1\leq j\leq k_i\}$ of $k_i$ variables each. It
has the following geometric interpretation: consider the system
\begin{equation}\label{syst}
F_i(t_1,\ldots,t_n):=\sum_{j=1}^{k_i}x_{ij}{\bf t}^{a_{ij}}=0  \ \ i=0,\ldots,n
\end{equation}
of Laurent polynomials in the variables $t_1, \dots, t_n.$
Here
${\bf t}^a$ stands for $t_1^{a_1}t_2^{a_2}\ldots t_n^{a_n}$ where
$a=(a_1,\ldots,a_n).$ The resultant $\RA$ vanishes on a particular
specialization of the $x_{ij}$ in an algebraically closed field
$K$ if the specialized system (\ref{syst}) has a common solution
in $\left(K\setminus\{0\}\right)^n$. See \cite{CLO,stu} for a
precise definition of $\RA$ and some basic facts.

Resultants are of fundamental importance for solving systems of polynomial equations and therefore have been extensively studied \cite{CLO,dan,EM,khe,min,stu}.
Recent research has focused on arithmetic aspects of this polynomial such as its {\it height} and its {\it Mahler measure} \cite{DH,KPS,som}.

Recall that the absolute height of $g:=\sum_\alpha c_{\alpha}X^\alpha\in\Cset [X_0,\ldots,X_n]$ is defined as $H(g):=\max\{|c_\alpha|,\,\alpha\in\Nset^k\},$ where
$k:=k_0+\ldots+k_n.$ Its (logarithmic) height is given by
\[h(g):=\log H(g)=\log\max\{|c_{\alpha}|,\,\alpha\in\Nset^k\}.\]
The Mahler measure of $g$ is defined as
\[ m(g):=\frac{1}{(2  \pi \ii)^{k}} \int_{\TT^{k}} \log |g(X_0, \dots, X_n) | \frac{ \dd X_0}{X_0} \dots \frac{ \dd X_n}{X_n},\]
where for $i=0,\ldots,n,\, \frac{\dd X_i}{X_i}$ is an abbreviation for
$\prod_{j=1}^{k_0}\frac{ \dd x_{ij}}{x_{ij}},$ and

\[\TT^{k} = \{(z_1,\dots,z_k) \in \Cset^{k} |\, |z_1| = \dots = |z_k| =1\}\] is the $k$-torus.

Some general relationships between the height and the Mahler measure are established in \cite[Chapter $3$]{EW} as well as \cite{KPS,som1}. In \cite{som}
upper bounds for both the height and the Mahler measure of resultants are presented.
However, very little seems to be known about the problem of explicitly computing  both the height and the Mahler measure of resultants. In the case of
heights, a first attempt was done in \cite{DH}, where the heights of resultants in low degree and one variable are calculated.

Jensen's formula gives a simple expression for the Mahler measure of a univariate polynomial as a function on its roots. However,
it is in general  a very hard problem to give an explicit closed formula for the Mahler measure of a multivariate polynomial.  The simplest examples are

\bigskip

\begin{theorem}
\begin{itemize}
\item \cite[Example 5]{smy1} \label{S1}
\begin{equation}
m(1+x+y) = \frac{3 \sqrt{3}}{4 \pi} \Lf(\chi_{-3},2)=  \Lf '
(\chi_{-3}, -1),
\end{equation}
where
\[{\Lf}  (\chi_{-3}, s) := \sum_{h=1}^\infty \frac{\chi_{-3}(h)}{h^s} \quad \mbox{with}\quad \chi_{-3}(h) := \left \{ \begin{array}{rlr} 1 & \mbox{if}\quad h \equiv1 \;& \mbox{mod} \; 3 \\ -1 & \mbox{if}\quad h \equiv -1\;&  \mbox{mod} \; 3 \\ 0 & \mbox{if}\quad h \equiv 0 \;& \mbox{mod} \; 3 \end{array} \right. \]
is the Dirichlet $\Lf$-series in the odd character of conductor $3$.
\item
Smyth also proved (see  \cite[Appendix 1]{Boy}):
\begin{equation}
m(1+x+y+z) = \frac{7}{2 \pi^2} \zeta(3),
\end{equation}
where $\zeta$ denotes the Riemann zeta function.
\end{itemize}
\end{theorem}

In this paper we focus in the explicit computation of the Mahler measure of $\RA$ in the case where the dimension of $N(\RA)$ (the Newton polytope of $\RA$) is small. We assume that the family of supports $\A_0,\ldots,\A_n$ is essential (see \cite[Sec.1]{stu}),
so $\RA$ is a polynomial of positive degree in the variables $x_{ij}.$ It is well-known (see \cite[Lemma 3.7]{EW})
that we always have $m(\RA)\geq0.$ The reason we focus in the dimension of the Newton polytope of the resultant and not in the number of variables and/or the size of the supports is due to some properties of the Mahler measure with respect to homogeneousness and changes of variables. For instance,  the Mahler measure of a homogeneous polynomial is the same as the Mahler measure of the corresponding dehomogenized polynomial. Moreover,

\begin{lemma}\cite[Lemma 7]{smy2} \label{varlemma}
Let $P({\bf y})$ be a $p$-variable polynomial, and let $V$ be a non-singular $p \times p$ integer matrix, then
\[m(P({\bf y})) = m(P({\bf y}^V)),\]
where ${\bf y}^V$ denotes $(\prod_j y_j^{v_{1j}}, \dots, \prod_j
y_j^{v_{pj}})$ for ${\bf y}=(y_1, \dots, y_p)$ and $V=\{ v_{ij}
\}$.
\end{lemma}

The whole situation may be summarized as follows: computing the Mahler measure of a polynomial whose Newton polytope has dimension $p$ is the same as computing the Mahler measure of a
$p$-variable polynomial. This is important because we may expect different kinds of formulas according to the number of variables (meaning the dimension of the Newton polytope). For speculations concerning this matter, see \cite{Lal}.

Evidence for this situation is Theorem \ref{mt} in Section \ref{22} which states that
$\RA$ has Mahler measure equal to zero if and only if its Newton polytope has dimension one, i.e. it is a segment.

This result shows that Mahler measures and heights behave differently in resultants. For instance, if we set
$n=1,\,\A_0=\{0,1\},\,\A_1=\{0,1,\ldots,\ell\},$ then it turns out that
\[{\rm Res}_{\A_0,\A_1}=\pm\sum_{j=0}^\ell(-1)^jx_{1j}x_{00}^{\ell-j}x_{01}^{j},\]
and from here it is easy to see that $h({\rm Res}_{\A_0,\A_1})=0.$ On the other hand, setting $y_j = (-1)^j x_{1j}x_{00}^{\ell-j}x_{01}^j$,
\[m({\rm Res}_{\A_0,\A_1})= m\left(\sum_{j=0}^\ell y_{j}\right).\]
The change of variables is allowable, because the $x_{1j}$ are algebraically independent,so we may apply Lemma \ref{varlemma}.

Dehomogeneizing, one obtains
\[m({\rm Res}_{\A_0,\A_1})=m(1+s_1+s_2+\ldots+s_{\ell}),\]
and this has been shown to be equal to $\frac12\log (\ell+1)-\frac{\gamma}{2}+O\left(\frac{\log (\ell+1)}{\ell+1}\right)$ as $\ell\to\infty,$
where $\gamma$ is the Euler--Mascheroni constant (see \cite{smy1}, and also \cite{R-VTV} for more estimates and generalizations).

Moreover, it is still unknown a characterization of all supports
$\A_0,\ldots,\A_n$ having $h(\RA)=0.$

In Section \ref{22} we deal with sparse resultants having Mahler measure zero.
Then we proceed to higher dimensions. In Section \ref{33} we focus on the case where the Newton polytope of the resultant has dimension two or three.
In Theorem \ref{dim2}, we compute the Mahler measure of resultants in dimension two, and in Theorem \ref{dim3},
we show that computing the Mahler measure of resultants in dimension three is essentially equivalent to the computation of Mahler measures of univariate trinomials. In Theorem \ref{tres} we compute the Mahler measure of trinomials having the same support.
In Section \ref{44} we compute the Mahler measure of a non trivial example in dimension four.

All the computations can be expressed  in terms of linear combinations of polylogarithms evaluated at algebraic numbers.
From the point of view of Mahler measure, it is natural to wonder why we would expect resultants to be a source of interesting examples of multivariate polynomials.

In \cite{D} Deninger established the relation between Mahler measure and regulators (see also \cite{RV,lal2}). More specifically, the Mahler measure of an irreducible polynomial $P \in \Qset[s_1,\dots ,s_p]$ is interpreted in terms of a special value of the regulator $\eta(s_1, \dots ,s_p)$ in $\mathcal{X}$, the projective variety determined by $\{ P = 0\}$. The regulator on the symbol $\{s_1,\dots, s_p\} \in K^M_p(\Cset(\mathcal{X}))\otimes \Qset$
is initially defined in the cohomology of $\mathcal{X} \setminus \{\mbox{poles and zeros of}\, s_i \}.$
A sufficient condition for extending it to the cohomology of $\mathcal{X}$ is that the tame symbols of the facets are trivial. In
that case the polynomial is called tempered (\cite{RV}).

If the symbol $\{s_1, \dots, s_p\}$ is trivial, then the tame symbols of the facets are trivial and $\eta(s_1, \dots ,s_p)$ is exact, and easily integrable by means of Stokes Theorem. This is the first step that may lead to a Mahler measure involving special values of polylogarithms (\cite{lal2}).

While the symbol is not necessarily trivial for a general polynomial, it is trivial for the case of sparse resultants. This is the content of Section \ref{55}. In Theorem \ref{delirio}, we show that resultants have trivial symbol, and so they are tempered polynomials. This fact suggests that the Mahler measure of resultants may be expressible in terms of combinations of polylogarithms and that we might expect results in the style of the ones from Sections \ref{33} and \ref{44} to be held in more generality.

\bigskip

\section{Resultants with Mahler measure equal to zero}\label{22}

The the main result of this section is the following:

\medskip

\begin{thm}\label{mt}
\[m (\RA )=0 \iff \dim(N(\RA))=1.\]
\end{thm}

\medskip

\begin{pf} Assume first that $\dim(N(\RA))=1.$ We use Proposition $4.1$ in \cite{stu}, which characterizes all families of essential supports $\A_0,\ldots,\A_n$ such that
the dimension of the Newton polytope is one: they must satisfy $k_i:=2,\,i=0,\ldots,n.$ It turns out that (\cite[Proposition 1.1]{stu}) in this case $\RA$ must be of the form
$\pm(X^{\lambda_1}-X^{\lambda_2}),$ with $\lambda_1,\lambda_2\in\Nset^k.$ It is very easy to see that polynomials of this kind have Mahler
measure zero because they are a monomial times the evaluation of $T-1$ in another Laurent monomial.

For the converse, assume that $m(\RA)=0,$. Recall that the resultant is a primitive polynomial in $\Zset[X_0,\ldots,X_n]$.  By Kronecker's Lemma  (see
for instance \cite[Theorem 3.10]{EW}), $\RA$ must be
a monomial times a product of cyclotomic polynomials evaluated in monomials. But $\RA$ is irreducible in $\Cset[X_0,\ldots,X_n]$ as it is the equation of an irreducible surface in the projective complex space (see \cite[Lemma 1.1]{stu}). Having its Mahler measure zero, the resultant must
be of the form $X^\alpha\pm X^\beta$ with $\alpha,\beta\in\Nset^k,$ i.e. a monomial times the polynomial $T\pm1$ evaluated at another Laurent monomial.
Hence, $\dim(N(\RA))=1.$
\end{pf}

\bigskip
\section{The Mahler measure of resultants in dimensions two and three}\label{33}
Now we would like to compute the Mahler measure of the systems having  $\dim N(\RA)>1.$
In order to do that, we first recall the following characterization of the dimension of the Newton polytope of the resultant:

\begin{theorem}\cite[Theorem $6.1$]{stu}\label{stt}
\[\dim(N(\RA))=k-2n-1,\]
where, as defined in the introduction, $k  = \sum_{i=0}^n k_i$.
\end{theorem}

We will compute the Mahler measure of the resultants having $\dim(N(\RA))=2.$ By the previous Theorem,
this property only holds in the case where there exists a unique $i_0$ such that $k_{i_0}=3$ and all other $k_i=2$, because the $k_i$ must be greater than 1 (see \cite[Theorem $1.1$]{stu}).

Suppose w.l.o.g. that $k_0=3$ and $k_1=k_2=\ldots=k_n=2.$ Consider any linear transformation in $SL(n,\Zset)$ which maps the directions in $\A_i$
to multiples $\eta_i{\bf e}_i$ of the unit vectors for $i=1,\ldots,n.$
After applying this transformation which does not change neither the Mahler measure nor the structure of $\RA,$ the original $F_i$'s defined in (\ref{syst}) look as follows:
\begin{equation}\label{2}
\begin{array}{ccl}
F_0(t_1,\ldots,t_n)&=&x_{01}{\bf t}^{a_{01}}+x_{02}{\bf t}^{a_{02}}+x_{03}{\bf t}^{a_{03}}, \\
F_1(t_1,\ldots,t_n)&=&x_{11}{t_1}^{\eta_1}-x_{12},\\
\ldots&\ldots&\ldots\\
F_n(t_1,\ldots,t_n)&=&x_{n1}{t_n}^{\eta_n}-x_{n2}.
\end{array}
\end{equation}

Let $\eta:=\eta_1+ \eta_2 + \ldots +\eta_n$.

\begin{thm}\label{dim2}
For systems having support as in (\ref{2}),
\[m(\RA)=  \eta \,\Lf '(\chi_{-3}, -1).\]
\end{thm}

\medskip
\begin{pf}
It is straightforward to verify that the resultant of (\ref{2}) is the following: for each $j=1,\ldots,n$ let $\xi_j$ run over the $\eta_j$-roots of unity. Then, it turns out that
$\RA$ equals, up to a monomial in the variables $x_{11},x_{21},\ldots,x_{n1},$
\[\prod_{j=1}^n\prod_{{\xi_j}^{\eta_j}=1}f_0\left(\xi_1\left(\frac{x_{12}}{x_{11}}\right)^{\frac{1}{\eta_1}},\xi_2\left(\frac{x_{22}}{x_{21}}\right)^{\frac{1}{\eta_2}},\ldots,
\xi_n\left(\frac{x_{n2}}{x_{n1}}\right)^{\frac{1}{\eta_n}}\right).\]
Let $V_i := \frac{x_{i2}}{x_{i1}}$. By Lemma \ref{varlemma}), 
\[M_j:= m\left(\prod_{\xi_j^{\eta_j}=1}f_0\left(\xi_1\left(\frac{x_{12}}{x_{11}}\right)^{\frac{1}{\eta_1}},\xi_2\left(\frac{x_{22}}{x_{21}}\right)^{\frac{1}{\eta_2}},\ldots,
\xi_n\left(\frac{x_{n2}}{x_{n1}}\right)^{\frac{1}{\eta_n}}\right)\right).\]
\[ =m\left(\prod_{\xi_j^{\eta_j}=1}f_0\left(\xi_1\left(V_1\right)^{\frac{1}{\eta_1}},\xi_2\left(V_2 \right)^{\frac{1}{\eta_2}},\ldots,
\xi_n\left(V_n \right)^{\frac{1}{\eta_n}}\right)\right).\]

Now, since $m(P(t_1,\ldots,t_n)) = m(P({t_1}^{p_1},\ldots,{t_n}^{p_n}))$ (by Lemma \ref{varlemma} onece again), we have
\[ M_j = m\left(\prod_{\xi_j^{\eta_j}=1}f_0\left(\xi_1 V_1 ,\xi_2 V_2,\ldots,\xi_n
V_n \right)\right).\]

Observe that coefficients of absolute value one can be absorbed by variables, so
$M_j = \eta_j\,m\left(f_0\left(V_1 ,V_2,\ldots,V_n \right)\right)$. Hence
\[M:=m(\RA)= \sum_{j=1}^n M_j = \eta\, m\left( f_0\left(V_1 ,V_2,\ldots,V_n\right)\right).\]

Now since $x_{01}$, $x_{02}$, and $x_{03}$ are algebraically independent, we may replace  $x_{01}V^{a_{01}}$, $x_{02}V^{a_{02}}$, and $x_{03}V^{a_{03}}$ by three independent variables $W_0$, $W_1$, and $W_2$,

\[ M = \eta\, m(W_0+W_1+W_2);\]
but this is just Smyth's result (Theorem \ref{S1}):

\[M = \eta\, m(1+x+y) = \eta \,\frac{3 \sqrt{3}}{4 \pi} \Lf(\chi_{-3},2)=
\eta\, \Lf '
(\chi_{-3}, -1).\]
\end{pf}

\bigskip
With the same proof as before, we can compute the Mahler measure of more general systems as follows. Consider an essential system of the form
\begin{equation}\label{gral}
\begin{array}{ccl}
F_0&=&x_{01}{\bf t}^{a_{01}}+x_{02}{\bf t}^{a_{02}}+\ldots+x_{0\ell}{\bf t}^{a_{0\ell}}, \\
F_1&=&x_{11}{t_1}^{\eta_1}-x_{12},\\
\ldots&\ldots&\ldots\\
F_n&=&x_{n1}{t_n}^{\eta_n}-x_{n2}.
\end{array}
\end{equation}

\begin{thm}\label{general}
With the notation established above, for systems as (\ref{gral}) we have
\[m(\RA)=\eta \,m(1+s_1+s_2+\ldots+s_{\ell-1}).\]
\end{thm}
As mentioned in the introduction, the Mahler measure of polynomials of the form $1+s_1+s_2+\ldots+s_p$ was estimated in \cite{smy1} and later in \cite{R-VTV}.

Now we would like to compute the Mahler measure of resultants having Newton polytope of dimension $3.$ According to Theorem \ref{stt}, we must consider
essentially the following two scenarios:
\begin{enumerate}
\item $k_0=4,\,k_1=k_2=\ldots=k_n=2.$ This is a system of the form (\ref{gral}), and hence we have that
\[m(\RA)=\eta\,m(1+s_1+s_2+s_3)=\eta\,\frac{7}{2\pi^2}\zeta(3),\] by Smyth's result (Theorem \ref{S1}).
\item $k_0=k_1=3,\,k_2=k_3=\ldots=k_n=2.$ This case is treated below.
\end{enumerate}

Let $k_0=k_1=3$ and $k_2=k_3=\ldots=k_n=2$. Consider a linear transformation in $SL(n,\Zset)$ which maps the directions in $\A_i$
to multiples $\eta_i{\bf e}_i$ of the unit vectors for $i=2,\ldots,n.$
\begin{equation}\label{3}
\begin{array}{ccl}
F_0&=&x_{01}{\bf t}^{a_{01}}+x_{02}{\bf t}^{a_{02}}+x_{03}{\bf t}^{a_{03}}, \\
F_1&=&x_{11}{\bf t}^{a_{11}}+x_{12}{\bf t}^{a_{12}}+x_{13}{\bf t}^{a_{13}}, \\
F_2&=&x_{21}{t_2}^{\eta_2}-x_{22},\\
\ldots&\ldots&\ldots\\
F_n&=&x_{n1}{t_n}^{\eta_n}-x_{n2}.
\end{array}
\end{equation}

Let $\alpha_{ij}\in\Zset$ be the first coordinate of the vector $a_{ij},\,
i=0,1,\, j=1,2,3,$ and
set as before $\eta:=\eta_2+\eta_3+\ldots+\eta_n.$ Consider the following system of supports:
\begin{equation}\label{je}
\A'_0:=\{\alpha_{01},\alpha_{02},\alpha_{03}\},\A'_1:=\{\alpha_{11},\alpha_{12},\alpha_{13}\}.
\end{equation}

Observe that the cardinalities of $\A'_0$ and $\A'_1$ must be at least two, otherwise the family $\A_0,\A_1,\ldots,\A_n$ would not be essential. We get the following.

\begin{thm}\label{dim3}
For systems like (\ref{3}) we have
\[m(\RA)=\eta\,m({\rm Res}_{\A'_0,\A'_1}).\]
\end{thm}

\begin{pf}
As in the proof of Theorem $6.2$ in \cite{stu}, it turns out that $\RA$ equals, up to a monomial factor, the product of the ${\rm Res}_{\A'_0,\A'_1}$
over all choices of roots of unity. We can then follow the same lines as in the proof of Theorem \ref{dim2} and conclude the claim.
\end{pf}

\bigskip
Therefore the computation of the Mahler measure of resultants in dimension three reduces to the computation of the Mahler measure of
univariate systems like (\ref{je}). Unfortunately, this does not seem to be very easy. In order to state our best result in that direction, we need to recall some facts about polylogarithms (see, for instance, \cite{Zag2}).

\begin{defn}
The {\em $q$th polylogarithm} is the function defined by the power series
\begin{equation}
 \Li_q(z) := \sum_{j=1}^\infty \frac{z^j}{j^q} \qquad z \in \Cset, \quad |z| <1.
\end{equation}
\end{defn}
This function can be continued analytically to $\Cset \setminus (1, \infty)$. Observe that $\Li_q(1)=\zeta(q)$ and $\Li_q(-1)=(2^{1-q}-1) \zeta(q)$.

In order to avoid discontinuities and to extend these functions to the whole complex plane, several modifications have been proposed. We will only need the cases $q=2, 3$. For $q=2$, we consider the Bloch--Wigner dilogarithm:
\begin{equation}
P_2(z) = D(z) := \im( \Li_2(z)) + \log |z| \arg(1-z).
\end{equation}
For $q=3,$ Zagier \cite{Zag2} proposes the following:
\begin{equation} \label{trilo}
 P_3(z) := \re\left( \Li_3(z) - \log |z| \Li_2(z) + \frac{1}{3} \log^2 |z| \Li_1(z) \right).
 \end{equation}

These functions are one-valued, real analytic in $\PP^1(\Cset) \setminus \{ 0, 1, \infty\}$, and continuous in $\PP^1(\Cset)$. Moreover, $P_q$ satisfies several functional equations, the simplest ones being, for $q=2$,
\begin{equation} \label{eq:prop1D}
D(\bar{z}) = -D(z), \qquad D(z) = -D(1-z) = -D \left( \frac{1}{z}\right),
\end{equation}
\begin{equation}\label{eq:prop3D}
 D(z) = \frac{1}{2} \left( D\left( \frac{z}{\bar{z}}\right)+ D\left( \frac{1-\frac{1}{z}}{1-\frac{1}{\bar{z}}}\right) + D\left( \frac{1-\bar{z}}{1-z}\right)\right).
\end{equation}

When $z$ has absolute value one, $D(z)$ has a particularly elegant expression:
\begin{equation} \label{eq:prop2D}
-2 \int_0^\theta \log|2 \sin t| \dd t =  D(\e^{2 \ii \theta}) = \sum_{j=1}^\infty \frac{\sin(2j\theta)}{j^2}.
\end{equation}
More about $D(z)$ can be found in \cite{Zag1}. For $q=3$, we have, for instance,
\begin{equation}
P_3(\bar{z}) =  P_3(z), \qquad P_3\left( \frac{1}{z} \right) =  P_3(z),
\end{equation}
\begin{equation}\label{eq:prop1P}
P_3(z) + P_3( 1 -z )+ P_3\left( 1 - \frac{1}{z} \right) = \zeta(3).
\end{equation}

We are now ready to state our result:

\begin{thm}\label{tres}
Suppose that $\A'_0=\A'_1,$ have both cardinality three.  W.l.o.g. we can suppose that $\A'_0=\{0,p,q\},$ with $p<q$ and $\gcd(p,q)=1.$

Then,
\[ m({\rm Res}_{\A'_0,\A'_1}) =\frac{2}{\pi^2}  \left(-p P_3(\varphi^q) -q P_3(-\varphi^p) +p P_3\left(\phi^q\right) + q P_3\left( \phi^p\right)\right)\]
where $\varphi$ is the real root of $x^{q}+x^{q-p}-1 = 0$ such that $0 \leq \varphi \leq 1,$ and $\phi$ is the real root of $x^{q} - x^{q-p} -1 =0$ such that $1 \leq \phi$.
\end{thm}

\medskip
\begin{pf}
All along this proof, we will write ${\rm Res}$ as short of ${\rm Res}_{\{0,p,q\},\{0,p,q\}}.$
First, we will show that
\[{\rm Res}(A+Bt^p+t^q,C+Et^p+t^q)= (C-A)^q-(EA-BC)^p(B-E)^{q-p}.\]

Let us set $f:=A+Bt^p+t^q$ and $g:=C+Et^p+t^q$. By using \cite[Ex $7$ Chapter 3]{CLO} we see that
\begin{equation}\label{cuk}
{\rm Res}(f,g)={\rm Res}(f,g-f)={\rm Res}(A+Bt^p+t^q,C-A+(E-B)t^p).\end{equation}
Let $\xi$ be a primitive $p$-th root of the unity, then all the roots of $C-A+(E-B)t^p$ are
$\xi^j\left(\frac{C-A}{B-E}\right)^{\frac{1}{p}},\ j=1,\ldots,p.$
By using the Poisson product formula for the computation of ${\rm Res}$ (see display (1.4) in Chapter $3$ of \cite{CLO}), we conclude that
(\ref{cuk}) equals
\begin{equation}\label{ck}
\begin{array}{c}
(-1)^{qp}(E-B)^q\prod_{j=1}^p f(\xi^j\left(\frac{C-A}{B-E}\right)^{\frac{1}{p}})\\
=
(-1)^{qp}(E-B)^q\prod_{j=1}^p\left(A+B\frac{C-A}{B-E}+\xi^{qj}\left(\frac{C-A}{B-E}\right)^{\frac{q}{p}}\right).
\end{array}
\end{equation}
The last product in (\ref{ck}) is of the form
\[\prod_{j=1}^p \alpha-\beta \xi^j=\alpha^p-\beta^p\]
with $\alpha=\frac{BC-AE}{B-E}$ and $\beta=-\left(\frac{C-A}{B-E}\right)^{\frac{q}{p}}$ (this is due to the fact that $\gcd(p,q)=1$).  So we get
that (\ref{ck}) equals
\[\begin{array}{l}
(-1)^{qp}(E-B)^q\left(\frac{(BC-AE)^p}{(B-E)^p}-(-1)^p\frac{(C-A)^q}{(B-E)^q}\right)\\
=(-1)^{qp}\left((-1)^q(BC-AE)^p(B-E)^{q-p}-(-1)^{p+q}(C-A)^q\right)\\
=(-1)^{qp+p+q+1}\left((C-A)^q-(AE-BC)^p(B-E)^{q-p}\right).
\end{array}
\]
The claim holds straightforwardly by noting that $qp+q+p+1=(q+1)(p+1)$ is even if $\gcd(p,q)=1.$

Now we have to compute the Mahler measure of \[(C-A)^q-(EA-BC)^p(B-E)^{q-p}.\]

After setting $C=C_1 A$, $E=E_1 B$ and dividing by $A^p$, we see that it is enough to consider the polynomial
$ A^{q-p}(C_1-1)^q - B^q (E_1-C_1)^p(1-E_1)^{q-p}.$

Now set $Z = A^{q-p}B^{-q}$ and divide by $B^q$. We need to compute
\begin{equation}\label{mll}
m(Z(C_1-1)^q -  (E_1-C_1)^p(1-E_1)^{q-p})
\end{equation}

By using Jensen's equality respect to the variable $Z$ and the fact that $m((C_1-1)^q)=0$, we deduce that (\ref{mll}) equals
\[ \frac{1}{(2 \pi \ii)^2} \int_{\TT^2} \log^+ \left | \frac{(E_1-1)^{q-p}(E_1-C_1)^p}{(C_1-1)^q} \right | \frac{\dd C_1}{C_1} \frac{\dd E_1}{E_1}, \]
where $\log^+|x|=\log|x|$ for $|x|\geq 1$ and zero otherwise.

Now write $C_1=E_1Y$. The expression above simplifies as follows:
\[ \frac{1}{(2 \pi \ii)^2} \int_{\TT^2} \log^+ \left | \frac{(E_1-1)^{q-p}(Y-1)^p}{(YE_1-1)^q} \right | \frac{\dd Y}{Y} \frac{\dd E_1}{E_1}. \]
Setting $Y = \e^{2 \ii \alpha}, E_1= \e^{2 \ii \beta}$, we have that this expression can be computed as follows:
\[ \frac{1}{\pi ^2} \int_{-\frac{\pi}{2}}^{\frac{\pi}{2}} \int_{-\frac{\pi}{2}}^{\frac{\pi}{2}} \log^+ \left | \frac{\sin^p{\alpha}\sin^{q-p}{\beta}}{\sin^q{(\alpha + \beta)}} \right | \dd \alpha \dd \beta\]
\begin{equation}\label{trigo}
 = \frac{2}{\pi ^2} \int_{-\frac{\pi}{2}}^{\frac{\pi}{2}} \int_0^{\frac{\pi}{2}} \log^+ \left | \frac{\sin^p{\alpha}\sin^{q-p}{\beta}}{\sin^q{(\alpha + \beta)}} \right | \dd \alpha \dd \beta.
\end{equation}
For $-\frac{\pi}{2} \leq \beta \leq 0$, set $\gamma = - \beta $. We can then
simplify (\ref{trigo}):
\[= \frac{2}{\pi^2} \int_0^{\frac{\pi}{2}} \int_0^{\frac{\pi}{2}} \log^+ \left | \frac{\sin^p{\alpha}\sin^{q-p}{\beta}}{\sin^q(\alpha + \beta)} \right | \dd \alpha \dd \beta \]
\[+ \frac{2}{\pi ^2} \int_0^{\frac{\pi}{2}} \int_0^{\frac{\pi}{2}} \log^+ \left | \frac{\sin^p{\alpha}\sin^{q-p}{\gamma}}{\sin^q(\alpha - \gamma)} \right | \dd \alpha \dd \gamma.
\]
Now we perform a change of variables. For the first term, write
\[ a = \frac{\sin \alpha}{\sin (\alpha+\beta)}, \quad b = \frac{\sin \beta}{\sin(\alpha+\beta)},\]
then
\[\dd \alpha \dd \beta = \frac{\dd a}{a} \frac{\dd b}{b}.\]
For the second term, set
\[ a = \frac{\sin \alpha}{\sin (\alpha - \gamma)}, \quad b = \frac{\sin \gamma}{\sin(\alpha - \gamma)},\]
then
\[\dd \alpha \dd \gamma = \frac{\dd a}{a} \frac{\dd b}{b}.\]

This change of variables has a geometric interpretation: we can think of $a$ and $b$ as the sides of a
triangle whose third side has length one. The side of length $a$ is opposite to the angle $\alpha$ and the side of length $b$ is opposite to $\beta$. This construction is possible because of the Sine Theorem.

\begin{figure}
\centering
\epsfig{file=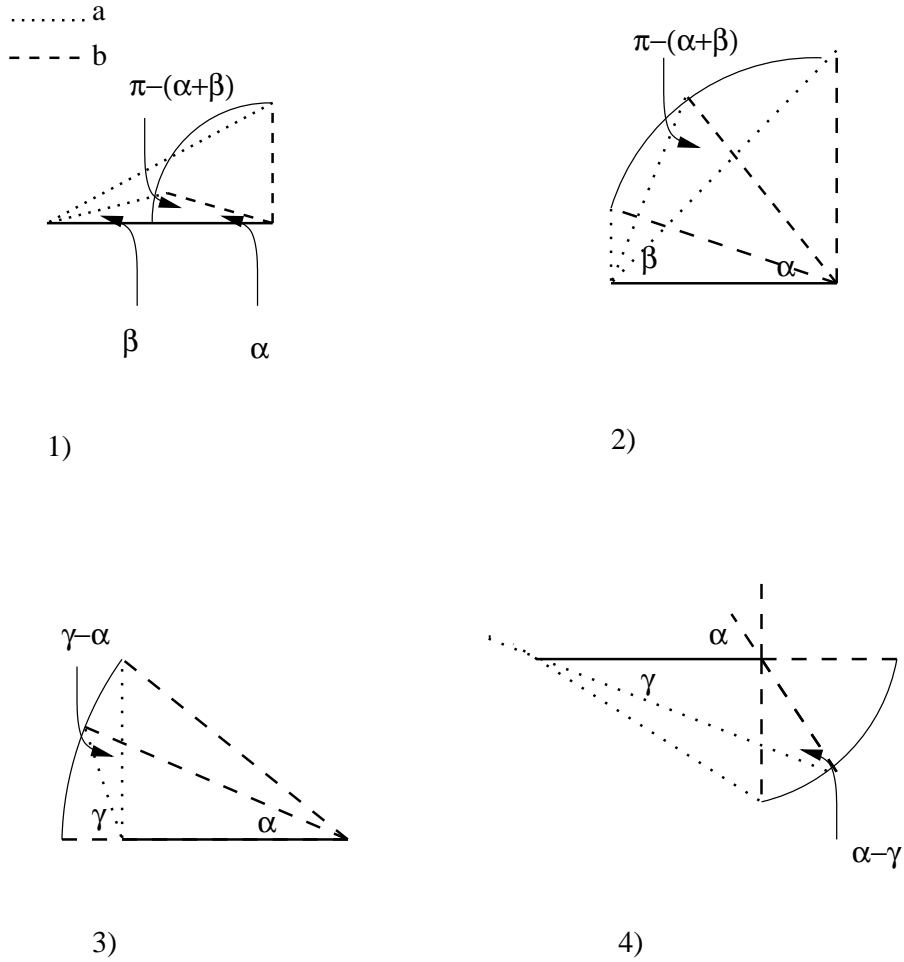}
\caption{1) Case when $0\leq b \leq 1$ in the first integral. 2) Case when $1 \leq b$ in the first integral.
3) Case when $\gamma \geq \alpha$ in the second integral. 4) Case when $\gamma \leq \alpha$ in the second integral.
  }
\label{dibuM}
\end{figure}

Figure \ref{dibuM} describes how the sides vary according to the angles.
The integral becomes the sum of four terms, each of them corresponding to each case in figure \ref{dibuM}.
\[ \frac{2}{\pi^2} \int_0^1 \int_{1-b}^{\sqrt{1+b^2}} \log^+ (a^p b^{q-p}) \frac{\dd a}{a} \frac{\dd b}{b}  +
\frac{2}{\pi^2} \int_1^\infty \int_{\sqrt{b^2-1}}^{\sqrt{1+b^2}} \log^+ (a^p b^{q-p}) \frac{\dd a}{a} \frac{\dd b}{b} \]
\[+ \frac{2}{\pi^2} \int_1^\infty \int_{b-1}^{\sqrt{b^2-1}} \log^+ (a^p b^{q-p}) \frac{\dd a}{a} \frac{\dd b}{b} +
\frac{2}{\pi^2} \int_0^\infty \int_{\sqrt{1+b^2}}^{1+b}\log^+ (a^p b^{q-p}) \frac{\dd a}{a} \frac{\dd b}{b} \]
\[ = \frac{2}{\pi^2} \int_0^1 \int_{1-b}^{1+b} \log^+ (a^p b^{q-p}) \frac{\dd a}{a} \frac{\dd b}{b}  + \frac{2}{\pi^2}\int_1^\infty \int_{b-1}^{1+b}\log^+ (a^p b^{q-p}) \frac{\dd a}{a} \frac{\dd b}{b}.\]

Now  write $c^{q-p}=a$ and $d^p = b$. Then the previous expression reduces to
\[
\frac{2p^2(q-p)^2}{\pi^2} \int_0^1 \int_{(1-d^p)^{\frac{1}{q-p}}}^{(1+d^p)^{\frac{1}{q-p}}} \log^+ (cd) \frac{\dd c}{c} \frac{\dd d}{d} + \frac{2p^2(q-p)^2}{\pi^2}\int_1^\infty \int_{(d^p-1)^{\frac{1}{q-p}}}^{(1+d^p)^{\frac{1}{q-p}}}\log^+ (cd) \frac{\dd c}{c} \frac{\dd d}{d}\]

\[ = \frac{2p^2(q-p)^2}{\pi^2} ( I_1 + I_2) \]

Let us compute $I_1$. Since the argument has the term $\log ^+ (cd)$,  we need to restrict the domain to the case $cd \geq 1$.
Now observe that since $0 \leq d \leq 1$, then $\frac{1}{d^{q-p}} \geq 1-d^p $. On the other hand, $ 1+d^p \geq \frac{1}{d^{q-p}} $ if and only if
\[ d^q + d^{q-p} - 1 \geq  0.\]
For future reference, let $\varphi$ be the unique root of $x^q+x^{q-p}-1 = 0$ in $[0,1]$.

\[ I_1 = \int_{\varphi}^1 \int_{\frac{1}{d}}^{(1+d^p)^{\frac{1}{q-p}}}\log (cd) \frac{\dd c}{c} \frac{\dd d}{d}
= \int_{\varphi}^1 \left. \frac{\log^2(cd)}{2}\right|_{\frac{1}{d}}^{(1+d^p)^{\frac{1}{q-p}}} \frac{\dd d}{d}\]
\[ = \int_\varphi^1 \frac{\log^2(d(1+d^p)^{\frac{1}{q-p}})}{2} \frac{\dd d}{d} \]

\begin{equation}\label{cuki}
= \int_{\varphi}^1 \left(\frac{\log^2 d}{2 d} + \frac{\log d\log (1+d^p)}{(q-p)d} + \frac{\log^2(1+d^p)}{2 (q-p)^2 d} \right ) \dd d.
\end{equation}
The first term in the integral (\ref{cuki}) is easy to integrate:
\[ \int_{\varphi}^1\frac{\log^2 d}{ d} \dd d = -\frac{\log^3 \varphi}{3}. \]
For the second term, we use the series expansion of $\log(1+x)$:
\[ \int_{\varphi}^1\frac{\log d\log (1+d^p)}{d} \dd d =  -\int_\varphi^1 \sum_{l=1}^\infty (-1)^l \frac{d^{pl-1}}{l} \log d \dd d\]
\[ = - \left . \sum_{l=1}^\infty (-1)^l \frac{d^{pl}}{pl^2} \log d \right |^1_\varphi + \int^1_\varphi \sum_{l=1}^\infty (-1)^l \frac{d^{pl-1}}{pl^2} \dd d \]
\[=  \frac{\log \varphi}{p} \Li_2(-\varphi^p) +\frac{1}{p^2}(\Li_3(-1) - \Li_3(-\varphi^p)).\]
We apply definition (\ref{trilo}) to conclude that this expression equals
\[ -\frac{1}{p^2}P_3(-\varphi^p) +\frac{(q-p)\log^3 \varphi }{3}  -\frac{3}{4p^2} \zeta(3).\]
Finally, we compute the third term of (\ref{cuki})
\[ \int_{\varphi}^1\frac{\log^2(1+d^p)}{ d} \dd d = \frac{1}{p} \int_{\frac{1}{2}}^{\varphi^{q-p}} \frac{\log^2 c}{c(1-c)} \dd c  \]
(setting $c = \frac{1}{1+d^p}$)
\[ =\frac{1}{p}\int_{\frac{1}{2}}^{\varphi^{q-p}} \log^2 c \left( \frac{1}{c} + \frac{1}{1-c} \right ) \dd c\]
\[ = \frac{(q-p)^3\log^3 \varphi}{3p} + \frac{\log^3 2}{3p} -\left. \frac{1}{p}\log^2 c \log(1-c)\right|^{\varphi^{q-p}}_{\frac{1}{2}}\]
\[+ \frac{1}{p}\int_{\frac{1}{2}}^{\varphi^{q-p}} \frac{2 \log c \log(1-c) }{c} \dd c \]
\[ = - \frac{(2q+p)(q-p)^2 \log^3 \varphi}{3p} - \frac{2 \log^3 2}{3p} - \frac{2}{p} \left . \sum_{l=1}^\infty \frac{c^{l}}{l^2} \log c \right |^{\varphi^{q-p}}_{\frac{1}{2}}\]
\[ + \frac{2}{p} \int_{\frac{1}{2}}^{\varphi^{q-p}}  \sum_{l=1}^\infty \frac{c^{l-1}}{l^2}\dd c \]
\[ = - \frac{(2q+p)(q-p)^2  \log^3 \varphi}{3p} - \frac{2 \log^3 2}{3p} - \frac{2(q-p)}{p}\log \varphi \Li_2(\varphi^{q-p})\]\[- \frac{2\log 2}{p} \Li_2\left( \frac{1}{2}\right) + \frac{2}{p}\Li_3(\varphi^{q-p}) - \frac{2}{p}\Li_3\left(\frac{1}{2}\right)\]
\[ =\frac{2}{p} P_3(\varphi^{q-p})- \frac{(q-p)^2 \log^3 \varphi}{3} -\frac{2}{p} P_3\left(\frac{1}{2}\right).\]

Then
\[I_1= \frac{1}{p(q-p)^2}P_3(\varphi^{q-p}) -\frac{1}{p^2(q-p)} P_3(-\varphi^p)- \frac{1}{p(q-p)^2}P_3\left(\frac{1}{2}\right)- \frac{3 \zeta(3)}{4 p^2 (q-p)}.\]

Let us compute $I_2$. As before, we need to restrict our domain to the case $cd \geq 1$. Since $1 \leq d$, we have $\frac{1}{d^{q-p}} \leq 1+d^p$ always. On the other hand, $d^p-1 \geq \frac{1}{d^{q-p}}$ if and only if
\[ d^q - d^{q-p} - 1 \geq  0.\]
For future reference, let $\phi$ be the unique root of $x^q-x^{q-p}-1 = 0$ in $[1, \infty)$.

Then
\[ I_2 = \int_1^{\phi} \int_{\frac{1}{d}}^{(1+d^p)^{\frac{1}{q-p}}}\log (cd) \frac{\dd c}{c} \frac{\dd d}{d} + \int_{\phi}^\infty \int_{(d^p-1)^{\frac{1}{q-p}}}^{(1+d^p)^{\frac{1}{q-p}}}\log (cd) \frac{\dd c}{c} \frac{\dd d}{d} = I_{21} + I_{22}.\]
We proceed to compute $I_{21}$,
\[  I_{21} = \int_1^{\phi}  \frac{\log^2(d(1+d^p)^{\frac{1}{q-p}})}{2} \frac{\dd d}{d} = \frac{1}{(q-p)^2}\int_{\phi^{-1}}^1 \frac{(\log(1+c^p) - q \log c )^2}{2} \frac{\dd c}{c}\]
(setting $c = \frac{1}{d} $)
\[ = \frac{1}{(q-p)^2}\int_{\phi^{-1}}^1 \left (  \frac{q^2\log^2 c}{2c} - \frac{q \log c \log (1+c^p)}{ c} + \frac{\log^2 (1+c^p)}{2 c}\right )  \dd c. \]

Using similar computations to those from $I_1$, we obtain
\[I_{21}= \frac{q^2 \log^3 \phi}{6(q-p)^2} + \frac{q \log^2 \phi \log \left( \frac{1+\phi^p}{\phi^p} \right)}{3(q-p)^2} + \frac{\log \phi \log^2 \left(\frac{1+ \phi^p}{\phi^p}\right)}{6(q-p)^2}+ \frac{q}{p^2(q-p)^2} P_3\left(-\frac{1}{\phi^p}\right)\]
\[ + \frac{1}{p(q-p)^2}P_3\left(\frac{\phi^p}{1+\phi^p}\right) - \frac{1}{p(q-p)^2}P_3\left(\frac{1}{2}\right) + \frac{3 q \zeta(3)}{4p^2(q-p)^2}.\]

For the case of $I_{22}$ we have:
\[ I_{22} = \int_{\phi}^\infty \frac{\log^2 (d(1+d^p)^{\frac{1}{q-p}}) - \log^2(d(d^p-1)^{\frac{1}{q-p}})}{2} \frac{\dd d}{d} \]
\[= \frac{1}{(q-p)^2} \int_0^{\phi^{-1}} \left( \frac{ \log^2(1+c^p)}{2c} - \frac{\log^2(1-c^p)}{2c} \right . \]
\begin{equation} \label{i22}
+ \left .\frac{q \log c \log(1-c^p)}{c} - \frac{q \log c \log(1+c^p)}{ c} \right) \dd c
\end{equation}
(setting $c= \frac{1}{d}$).

Now we compute each of the terms in equation  (\ref{i22}):
\[\int_0^{\phi^{-1}} \frac{ \log^2(1+c^p)}{c} \dd c=  -\frac{\log^2 \left( \frac{\phi^p}{1+\phi^p}\right) \log\phi}{3} - \frac{2}{p}P_3\left(\frac{\phi^p}{1+\phi^p}\right) +\frac{2}{p} \zeta(3).\]

\[\int_0^{\phi^{-1}} \frac{ \log^2(1-c^p)}{c} \dd c= \frac{1}{p}\int_{\phi^{-q}}^1 \frac{\log^2 f}{1-f} \dd f  \]
(setting $f = 1-c^p$)
\[ = - \frac{1}{p}\left .\log^2 f \log(1-f) \right |^1_{\phi^{-q}} + \frac{1}{p} \int_{\phi^{-q}}^1 \frac{2 \log f \log(1-f)}{f} \dd f
\]
\[ = -q^2 \log^3 \phi -\frac{2 q}{p} \log \phi \Li_2\left(\frac{1}{\phi^q}\right) + \frac{2}{p} \zeta(3) - \frac{2}{p} \Li_3\left(\frac{1}{\phi^q}\right) \]
\[=-\frac{q^2 \log^3 \phi}{3} - \frac{2}{p} P_3\left(\frac{1}{\phi^q}\right)+ \frac{2}{p} \zeta(3).\]
\[\int_0^{\phi^{-1}} \frac{\log c \,\log(1-c^p)}{c} \dd c =\frac{\log \phi}{p} \Li_2\left(\frac{1}{\phi^p}\right) + \frac{1}{p^2}\Li_3\left(\frac{1}{\phi^p}\right) \]
\[= \frac{1}{p^2}P_3\left(\frac{1}{\phi^p}\right) - \frac{q \log^3 \phi}{3}.\]

\[\int_0^{\phi^{-1}} \frac{\log c \,\log(1+c^p)}{c} \dd c =  \frac{\log \phi}{p} \Li_2\left(-\frac{1}{\phi^p}\right) + \frac{1}{p^2}\Li_3\left(-\frac{1}{\phi^p}\right) \] \[=
\frac{1}{p^2}P_3\left( - \frac{1}{\phi^p}\right) + \frac{\log^2\phi \log\left(\frac{1+\phi^p}{\phi^p}\right)}{3}.\]

Putting all the terms together,
\[ (q-p)^2 I_{22} = -\frac{\log^2 \left( \frac{\phi^p}{1+\phi^p}\right) \log\phi}{6} - \frac{q^2 \log^3 \phi}{6} - \frac{q \log^2 \phi \log \left(\frac{1+\phi^p}{\phi^p}\right)}{3}\]
\[- \frac{1}{p}P_3\left(\frac{\phi^p}{1+\phi^p}\right) + \frac{1}{p}P_3\left( \frac{1}{\phi^q}\right) + \frac{q}{p^2} P_3\left( \frac{1}{\phi^p}\right) - \frac{q}{p^2} P_3\left(- \frac{1}{\phi^p}\right) \]
and hence,
\[(q-p)^2I_2 =  - \frac{1}{p}P_3\left(\frac{1}{2}\right) + \frac{3 q \zeta(3)}{4p^2} + \frac{1}{p}P_3\left( \frac{1}{\phi^q}\right) + \frac{q}{p^2} P_3\left( \frac{1}{\phi^p}\right).  \]

Now we can conclude:
\[ I_1 + I_2 = \frac{1}{p(q-p)^2}P_3(\varphi^{q-p}) -\frac{1}{p^2(q-p)} P_3(-\varphi^p)\]
\[ - \frac{2}{p(q-p)^2} P_3\left (\frac{1}{2} \right) + \frac{3 \zeta(3)}{4p(q-p)^2} + \frac{1}{p(q-p)^2}P_3\left(\phi^q\right) + \frac{q}{p^2(q-p)^2} P_3\left( \phi^p\right).  \]

Let us note that $ 2 P_3\left (\frac{1}{2} \right) + P_3(-1) =  \zeta(3)$ because of equation (\ref{eq:prop1P}), hence,
\[ P_3\left (\frac{1}{2} \right) = \frac{7}{8} \zeta(3). \]

Then

\[I_1+I_2 =\frac{1}{p(q-p)^2}P_3(\varphi^{q-p}) -\frac{1}{p^2(q-p)} P_3(-\varphi^p)\]
\[ - \frac{ \zeta(3)}{p(q-p)^2} + \frac{1}{p(q-p)^2}P_3\left(\phi^q\right) + \frac{q}{p^2(q-p)^2} P_3\left( \phi^p\right).  \]

Now, we use equation (\ref{eq:prop1P}) again in order to obtain $P_3(\varphi^q) + P_3(1-\varphi^q) + P_3(1-\varphi^{-q}) = \zeta(3)$ from where
\[P_3(\varphi^q) + P_3(\varphi^{q-p}) + P_3(-\varphi^{-p}) = \zeta(3),\]
so
\[ P_3(\varphi^{q-p}) = \zeta(3) - P_3(-\varphi^p) - P_3(\varphi^q). \]

Hence, we have
\[ I_1+I_2 = -\frac{1}{p(q-p)^2}  P_3(\varphi^q) -\frac{q}{p^2(q-p)^2} P_3(-\varphi^p)\]
\[  + \frac{1}{p(q-p)^2}P_3\left(\phi^q\right) + \frac{q}{p^2(q-p)^2} P_3\left( \phi^p\right),  \]
which proves our claim.
\end{pf}

\bigskip
\section{An example in dimension $4$}\label{44}
We would like to study one more example, which is a particular case of a 4-dimensional resultant. Let us set $n=2$ and
$$\A_0=\A_1=\A_2=\A:=\{(0,0),(1,0),(0,1)\}.$$
We will use a formula due to Cassaigne and Maillot

\begin{theorem} \cite[Proposition 7.3.1]{Mai},
\begin{equation}\label{eq:Mai}
\pi m(a+bx+cy) = \left \{ \begin{array}{lr} D\left(\left|\frac{a}{b}\right| \e^{\ii \gamma}\right) + \alpha \log |a| + \beta \log |b| +\gamma \log |c|   & \mbox{if}\, \triangle  \\\\ \pi \log \max \{ |a|, |b|, |c|\} &  \mbox{if not}\, \triangle \end{array} \right.
\end{equation}
Where $\triangle$ stands for the statement that $|a|$, $|b|$, and
$|c|$ are the lengths of the sides of a triangle; and $\alpha$,
$\beta$, and $\gamma$ are the angles that are opposite to the sides of lengths
$|a|$, $|b|$ and $|c|$ respectively.
\end{theorem}

\begin{thm}
\[ m( \Res_{\A,\A,\A})  = \frac{9 \zeta(3)}{2 \pi^2}.\]
\end{thm}

\begin{pf}
In order to simplify the notation, we will use the variables $a,b,c,\dots$ instead of the $x_{ij}$'s.
\[\Res_{\A,\A,\A} = \det \left( \begin{array}{ccc} a & b & c \\ d & e & f \\ g & h & i \end{array} \right ).\]

Now, let us proceed to eliminate homogeneous variables:
\[m\left( \left| \begin{array}{ccc} a & b & c \\ d & e & f \\ g & h & i \end{array} \right | \right) =
m\left( \left| \begin{array}{ccc} 1 & 1 & 1 \\ d & e & f \\ g & h & i \end{array} \right | \right) =
m\left( \left| \begin{array}{ccc} 1 & 1 & 1 \\ 1 & e & f \\ 1 & h & i \end{array} \right | \right) \]
\[=
m\left( \left| \begin{array}{ccc} 1 & 0 & 0 \\ 1 & e-1 & f-1 \\ 1 & h-1 & i-1 \end{array} \right | \right) = m((e-1)(i-1)-(f-1)(h-1)). \]

Let us observe that
\begin{equation}\label{inter}
m( (x-1)(y-1) - (z-1)(w-1)) = m((x-1) y +(1-z) w +(z-x)).
\end{equation}
Hence we can think of the polynomial $(x-1) y +(1-z) w +(z-x)$ as a linear polynomial in the variables $y$, $w$, whose coefficients are in $\Zset[x,z]$. Because of the iterative nature of the definition of Mahler measure, we can choose to integrate first respect to the variables $y$ and $w$, regarding $x$ and $z$ as parameters. If we do that, we obtain formula (\ref{eq:Mai}) with the sides of the triangle equal to $|x-1|$, $|z-1|$, and $|z-x|$.

Now in order to compute the Mahler measure, we still need to integrate this formula respect to $x$ and $z$. Set $x = \e^{ 2\ii \alpha}$ and $z = \e^{ 2\ii \beta}$.
This notation is consistent with the names for the angles of the triangle because of the Sine Theorem (see Figure \ref{carlos1}).

\begin{figure}
\centering
\epsfig{file=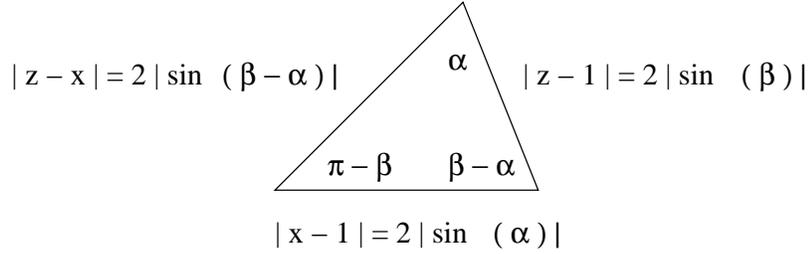}
\caption{We always obtain a triangle for this case.}
\label{carlos1}
\end{figure}

We obtain that (\ref{inter}) is equal to
\[\frac{2}{\pi^3} \int_0^\pi \int_0^\beta D\left(\frac{\sin \beta}{\sin \alpha} \e^{\ii ( \beta-\alpha)}\right)
+ \alpha \log|2 \sin \alpha| + (\beta - \alpha) \log | 2 \sin (\beta -\alpha)|\]
\[ + (\pi -\beta) \log |2 \sin \beta| \dd \alpha \dd \beta.\]
First, we integrate the terms involving logarithms:
\[\int_0^\pi \int_\alpha^\pi \alpha \log|2 \sin \alpha| \dd \beta \dd \alpha = \int_0^\pi \alpha(\pi-\alpha) \log|2 \sin \alpha| \dd \alpha\]
\[ =  \int_0^\pi \alpha(\pi-\alpha) \log\left|1-\e^{2\ii\alpha}\right| \dd \alpha =  -\int_0^\pi \alpha(\pi-\alpha) \re \sum_{j=1}^\infty\frac{\e^{2j\ii\alpha}}{j} \dd \alpha.\]
Because of
\[\int_0^\pi (\pi\alpha-\alpha^2) \frac{\cos(2j\alpha)}{j} \dd \alpha= \left. (\pi \alpha - \alpha^2) \frac{\sin(2j \alpha)}{2j^2}\right |^\pi_0 - \int_0^\pi (\pi  - 2\alpha) \frac{\sin(2j \alpha)}{2j^2} \dd \alpha\]
\[=  \left. (\pi - 2 \alpha) \frac{\cos(2j \alpha)}{4j^3}\right |^\pi_0 -\int_0^\pi (- 2) \frac{\cos(2j \alpha)}{4j^3} \dd \alpha = - \frac{\pi}{2 j^3}, \]
we conclude that
\[ \int_0^\pi \alpha(\pi-\alpha) \log|2 \sin \alpha| \dd \alpha= \frac{\pi \zeta(3)}{2}.\]

Then
\[\int_0^\pi \int_0^\beta(\pi -\beta) \log |2 \sin \beta| \dd \alpha \dd \beta= \int_0^\pi \beta(\pi - \beta) \log |2 \sin \beta| \dd \beta =\frac{\pi \zeta(3)}{2},\]
by analogy with the case of $\alpha$.

By setting $\gamma=\beta-\alpha$ in the third logarithmic term,
\[\int_0^\pi \int_0^\beta(\beta - \alpha) \log | 2 \sin (\beta -\alpha)|\dd \alpha \dd \beta =
\int_0^\pi \int_0^{\pi- \gamma} \gamma \log | 2 \sin \gamma|\dd \alpha \dd \gamma \]
\[ =\int_0^\pi \gamma(\pi - \gamma) \log |2 \sin \gamma| \dd \gamma = \frac{\pi \zeta(3)}{2}.\]

On the other hand, we need to evaluate
\begin{equation}\label{intdilog}
\int_0^\pi \int_0^\beta D\left(\frac{\sin \beta}{\sin \alpha} \e^{\ii ( \beta-\alpha)}\right)\dd \alpha \dd \beta.
\end{equation}

Using equation (\ref{eq:prop3D}),
\[D\left(\frac{\sin \beta}{\sin \alpha} \e^{\ii ( \beta-\alpha)}\right) = D\left(\frac{1-z}{1-x}\right) = \frac{1}{2} \left( D\left(\frac{z}{x} \right)+ D(x) + D\left( \frac{1}{z}\right)\right) \]
\[=\frac{1}{2} \left( D(\e^{2 \ii (\beta-\alpha)})+ D(\e^{2 \ii \alpha}) + D(\e^{-2 \ii \beta}) \right).\]

Then integral (\ref{intdilog}) is the sum of three terms. We proceed to compute each of them:
\[\int_0^\pi \int_\alpha^\pi D(\e^{2 \ii \alpha}) \dd \beta \dd \alpha = \int_0^\pi (\pi - \alpha) D(\e^{2 \ii \alpha}) \dd \alpha \]
\[=\int_0^\pi (\pi - \alpha) \sum_{j=1}^\infty \frac{\sin ( 2 j \alpha)}{j^2} \dd \alpha.\]

But
\[\int_0^\pi (\pi - \alpha)\frac{\sin ( 2 j \alpha)}{j^2} \dd \alpha = \left. -(\pi  - \alpha) \frac{\cos(2j \alpha)}{2j^3}\right |^\pi_0 - \int_0^\pi \frac{\cos(2j \alpha)}{2j^3} \dd \alpha = \frac{\pi}{2 j^3}.\]

Then
\[ \int_0^\pi (\pi - \alpha) D(\e^{2 \ii \alpha}) \dd \alpha = \frac{\pi \zeta(3)}{2}.\]

The other terms can be computed in a similar fashion:
\[\int_0^\pi \int_0^\beta D(\e^{-2 \ii \beta}) \dd \alpha \dd \beta = - \int_0^\pi \beta D(\e^{2 \ii \beta})\dd \beta = \frac{\pi \zeta(3)}{2},\]
\[\int_0^\pi \int_0^\beta D(\e^{2 \ii (\beta-\alpha)}) \dd \alpha \dd \beta= \int_0^\pi \int_0^{\pi- \gamma} D(\e^{2 \ii \gamma }) \dd \alpha \dd \gamma\]
\[ = \int_0^\pi (\pi - \gamma)D(\e^{2 \ii \gamma }) \dd \gamma = \frac{\pi \zeta(3)}{2}.\]

Thus, we conclude that
\[ m( (x-1)(y-1) - (z-1)(w-1)) = \frac{2}{\pi^3} \left( \frac{3}{2} \frac{\pi \zeta(3)}{2} + 3  \frac{\pi \zeta(3)}{2} \right) = \frac{9 \zeta(3)}{2 \pi^2}.\]
\end{pf}

\bigskip
\section{Resultants are tempered polynomials}\label{55}
In this section we will leave the elementary approach given above and turn instead to study algebraic properties of resultants and Mahler measures in the context of Milnor $K$-theory. We will sketch the relation here and we refer to \cite{D}, \cite{RV}, and \cite{lal2} for precise details.

For $P \in \Cset[s_1, \dots, s_p]$ irreducible, we write
\[P(s_1,\dots,s_p)= \sum_{i\geq 0} a_i(s_1,\dots,s_{p-1}) s_p^i.\]
Let
\[ P^*(s_1,\dots,s_{p-1}) = a_{i_0}(s_1,\dots,s_{p-1}),\]
the main non-zero coefficient respect to $s_p$. Let ${\mathcal X}$ be the zero set $\{P(s_1,\dots, s_p) =0 \}$.

By applying Jensen's formula to the Mahler measure of $P$ respect to the variable $s_p$, it is possible to write
\[ m(P)=m(P^*) - \frac{1}{(2\ii \pi)^{p-1}} \int_\Gamma \log|s_p|  \frac{\dd s_1}{s_1} \dots  \frac{\dd s_{p-1}}{s_{p-1}} \]
where $\Gamma = \{P(s_1,\dots, s_p) =0 \} \cap \{|s_1|=\dots=|s_{p-1}|=1, |s_p| \geq 1\}$.

From this formula, Deninger \cite{D} establishes
\begin{equation}
m(P)=m(P^*) - \frac{1}{(2\ii \pi)^{p-1}} \int_\Gamma \eta(s_1,\dots, s_p)
\end{equation}
where $\eta(s_1,\dots, s_p)$ is certain $\Rset(p-1)$-valued smooth $p-1$-form in $\mathcal{X}(\Cset) \setminus S$. Here, $S$ denotes the set of poles and zeros of the functions $s_1, \dots, s_p$.

For example, in two variables, $\eta$ has the following shape:
\[ \eta(x,y) = \log|y| \ii\dd \arg x - \log|x| \ii\dd \arg y.\]
$\eta$ is a closed form that is multiplicative and antisymmetric in the variables $s_1, \dots s_p$. Therefore, it is natural to think of $\eta$ as a function on $\bigwedge^p\Cset(\mathcal{X})\otimes \Qset$ (we tensorize by $\Qset$ because $\eta$ is trivial in torsion elements). Even more, $\eta(1-s, s, s_3, \dots, s_p)$ is an exact form. We may also see $\eta(s_1, \dots, s_p)$ as a class in the (DeRham) cohomology of $\mathcal{X} \setminus S$. This situation allows us to think of the cohomological class of $\eta$ as a function in the Milnor $K$-theory group $K^M_p(\Cset(\mathcal{X})) \otimes \Qset$. Recall that for a field $F$ the Milnor $K$-theory group is given by
\[ K^M_p(F) :=\bigwedge^p F^*/\left< (1-s_1)\wedge s_1 \wedge \dots \wedge s_p, s_i \in F^* \right >\]

If we can extend this class to the cohomology of $\mathcal{X}$, $\eta$ becomes a regulator. In certain cases, seeing the Mahler measure as a regulator allows us to explain its relation to special values of L-functions via Beilinson's conjectures and similar results.

The condition that the class of  $\eta(s_1, \dots, s_p)$ be extended to $\mathcal{X}$ is given by the triviality of the tame symbols in the Milnor $K$-theory. A stronger condition is that $\eta(s_1, \dots, s_p)$ is exact. Since $\eta$ is defined in $K^M_p(\Cset(\mathcal{X}))\otimes \Qset$, $\eta(s_1, \dots, s_p)$ is exact if the symbol $\{s_1, \dots, s_p\}$ is trivial in $K^M_p(\Cset(\mathcal{X}))\otimes{\Qset}$.

This is a very special condition that is not true for a general polynomial. However, it is true for resultants:

\begin{thm}\label{delirio} The symbol
\begin{equation}\label{symbol}
\{x_{01}, \dots, x_{0k_0},\dots,x_{n1}, \dots, x_{nk_n}\} \in K^M_k(\Cset(\mathcal{X}))\otimes{\Qset}
\end{equation}
is trivial.
\end{thm}


In order to prove Theorem \ref{delirio}, we will need the following
\begin{lem}\label{Lemma2}
Consider a $p$-variable polynomial \[ P(u_1,\dots, u_p) := \sum x_{\bf i} u^{\bf i}.\]
Here {\bf i} is a multiindex, $u^{\bf i} = u_1^{i_1} \dots u_p^{i_p}$.

Let $E$ be a field containing the $x_{\bf i}'s$ and let $\alpha_1, \dots \alpha_p$ in $E$. Then $\bigwedge x_{\bf i}$ is of the form
 \[P(\alpha_1, \dots, \alpha_p)  \wedge  \bigwedge a_i + \sum_{i,h}r_{i,h} \alpha_i \wedge \bigwedge c_{i,h}   + \sum_{j,l}s_{j,l} b_j\wedge(1-b_j) \wedge \bigwedge b'_{j,l},\]
where $a_i, b_j, b'_{j,l},c_{i,h}$ are elements of $E$ and $r_{i,h},s_{j,l}$ are integral numbers.
\end{lem}

\begin{pf} First we prove the case for which $p=1$. Let us write
\[ P(\alpha) = x_0 \left ( 1 + \frac{x_1}{x_0} \alpha  \left ( 1 + \frac{x_2}{x_1} \alpha  \dots \left(1+\frac{x_q}{x_{q-1}}\alpha \right)\right) \right ).\]

After setting $y_0=x_0$ and $y_i = \frac{x_i}{x_{i-1}}$ for $i >0$, we obtain
\[ P(\alpha) = y_0( 1 + y_1 \alpha (1 + y_2 \alpha \dots (1+y_q \alpha ))),\]
and
\[ x_0 \wedge \dots \wedge x_q = y_0 \wedge \dots \wedge y_q.\]

We may introduce $\alpha$ in  the last place of the wedge product of $y_i$:
\[y_0 \wedge \dots \wedge y_q = y_0 \wedge \dots \wedge y_{q-1}\wedge (y_q \alpha) -y_0 \wedge \dots \wedge y_{q-1}\wedge  \alpha.\]
Now we introduce $\alpha$ in the second to the last place
\[ = y_0 \wedge \dots \wedge y_{q-1}(1+y_q \alpha) \wedge (y_q \alpha) -y_0 \wedge \dots \wedge (1+y_q \alpha) \wedge (y_q \alpha) \]
\[ -y_0 \wedge \dots \wedge y_{q-1}\wedge  \alpha \]

\[ = y_0 \wedge \dots \wedge y_{q-2} \wedge y_{q-1}\alpha (1+y_q \alpha) \wedge (y_q \alpha)  - y_0 \wedge \dots \wedge y_{q-2} \wedge \alpha \wedge y_q   \]
\[-y_0 \wedge \dots \wedge (1+y_q \alpha) \wedge (y_q \alpha)  -y_0 \wedge \dots \wedge y_{q-1}\wedge  \alpha.\]

Then we introduce  $\alpha$ in the third to the last place. We continue in this fashion
\[ = y_0 \wedge \dots \wedge y_{q-2}(1+y_{q-1}\alpha (1+y_q \alpha) )  \wedge y_{q-1}\alpha (1+y_q \alpha) \wedge (y_q \alpha)\]
\[ - y_0 \wedge \dots \wedge(1+y_{q-1}\alpha (1+y_q \alpha) )  \wedge y_{q-1}\alpha (1+y_q \alpha) \wedge (y_q \alpha)  \]
\[- y_0 \wedge \dots \wedge y_{q-2} \wedge \alpha \wedge y_q   -y_0 \wedge \dots \wedge (1+y_q \alpha) \wedge (y_q \alpha)  \]
\[ -y_0 \wedge \dots \wedge y_{q-1}\wedge  \alpha .\]

After $q-1$ more steps, we get
\[ P(\alpha)  \wedge \bigwedge a_i + \sum_{h=1}^q \pm ( b_h \wedge (1-b_h) \wedge \bigwedge b'_{h,l} ) \]
\[ - \sum_{h=1}^q y_0 \wedge \dots \wedge y_{h-1} \wedge \alpha \wedge y_{h+1} \wedge \dots
\wedge y_q,\]
which proves the claim in this case.

 Let us now consider $p>1$. We use induction on $p$. Suppose that the claim is true for $p-1$. Then we may write
\[ P(\alpha_1,\dots, \alpha_p) = \sum P_i(\alpha_1, \dots, \alpha_{p-1}) \alpha_p^i.\]
By the inductive hypothesis we obtain
\[ \bigwedge x_{\bf i} = P_0(\alpha_1, \dots, \alpha_{p-1}) \wedge \dots \wedge P_h(\alpha_1, \dots, \alpha_{p-1})  \wedge \bigwedge a_i  \]
 \[+ \sum_{i,h}r_{i,h} \alpha_i \wedge \bigwedge c_{i,h}   + \sum_{j,l}s_{j,l} b_j\wedge(1-b_j) \wedge \bigwedge b'_{j,l}.\]
Now apply the case $p=1$ to the first term in order to obtain the desired result.

In the notation of this proof we used that the coefficients were all different from zero. The proof is easily adapted to the case when some coefficients are zero. We are going to use this lemma in full generality.
\end{pf}

\bigskip
\begin{pf*}{PROOF (Theorem \ref{delirio})} By definition of Milnor's $K$-theory, it is sufficient to work in $\bigwedge^{k} \Cset(\mathcal{X})^*\otimes \Qset$.  Consider the variety
\[\mathcal{Y}=\{z-\RA=0\}\subset\Cset^{k+1}.\]
Let $H:=  \Cset(\mathcal{Y})$ and $E$ be a finite extension of $H$ containing all the  roots that are common to the last $n$ polynomials $F_1, \dots, F_n$ in an algebraically closed field containing $H.$

We will prove that the symbol determined by
\[ z \wedge \bigwedge x_{0j} \wedge \bigwedge x_{1j} \wedge \dots \wedge \bigwedge x_{nj}  \]
is trivial in $\bigwedge^{k+1}E^* \otimes \Qset$. This fact implies that the corresponding symbol in $K^M_{k+1}(E)\otimes \Qset$ is trivial.
Now $j:H\hookrightarrow E$ induces $j:K^M_*(H) \rightarrow K^M_*(E)$ whose kernel is finite (see \cite{BT}). Then  $K^M_*(H)\otimes \Qset \hookrightarrow K^M_*(E)\otimes \Qset$ is injective and we conclude that the symbol must be trivial in $\bigwedge^{k+1}H^*\otimes \Qset$.

The triviality of this symbol implies that
\[\bigwedge x_{0j} \wedge \bigwedge x_{1j} \wedge \dots \wedge \bigwedge x_{nj}\] is trivial in $\bigwedge^{s} \Cset(\mathcal{X})^*\otimes \Qset$, because it is the image by the tame symbol morphism respect to the valuation determined by $z=0$, \cite{Mil} (in the language of Newton polytopes, $\RA$ corresponds to a facet of $z-\RA$).

We will use the Poisson product formula for sparse resultants (Theorem $1.1$ in \cite{PS} and its refinement in \cite[Theorem 8]{min}):
\begin{equation}\label{pescado}
\RA = \prod_{\eta}{\rm Res_\eta}^{d_\eta}\prod_{\tilde{\alpha}} F_0(\tilde{\alpha}),
\end{equation}
where $\eta$ runs over all the maximal facets of the Newton polytope associated to the Minkowski sum $\A_1+\ldots+\A_n,$ $\rm Res_\eta$ stands for the \textit{facet resultant} associated to
(\ref{syst}) and the facet $\eta,\, d_\eta$ is a non negative integer number, and $\tilde{\alpha}$ runs over all the common solutions of the system
$F_1=\ldots=F_n=0$ in $\left(E\setminus\{0\}\right)^n.$

Let us proceed by induction on $n$. For $n=1$ we have, $F_0=x_{01}t_1^{a_{01}} + \dots + x_{0k_0}t_1^{a_{0k_0}}$ and $F_1=x_{11}t_1^{a_{11}} + \dots + x_{1k_1}t_1^{a_{1k_1}}$, where
we assume w.l.o.g. that $a_{01}<a_{02}<\dots <a_{0k_0}$ and $a_{11}<a_{12}<\dots <a_{1k_1}.$ Also, we can suppose w.l.o.g. that $a_{01}=a_{11}=0.$ This is due to the fact that the resultant is invariant under translations
of the $\A_i$'s.
Then,
\[\Res(F_0, F_1) = x_{1k_1}^{d}\prod F_0(\alpha), \]
where $\alpha$ runs over the roots of $F_1$ in $E\setminus\{0\},$ and $d$ is a positive integer.
Let $\tilde{F}_1:= x_{12}t_1^{a_{12}} + \dots + x_{1k_1}t_1^{a_{1k_1}}.$
Then
\[ z\wedge \bigwedge x_{0j} \wedge \bigwedge x_{1j} = \sum F_0(\alpha) \wedge \bigwedge_{j=1}^{k_0} x_{0j} \wedge \tilde{F}_1(\alpha) \wedge \bigwedge_{j=2}^{k_1} x_{1j}\]
\[ + d \, x_{1k_1}\wedge \bigwedge_{j=1}^{k_0} x_{0j} \wedge \tilde{F}_1(\alpha) \wedge \bigwedge_{j=2}^{k_1} x_{1j}.\]
The second term is zero because it contains two copies of $x_{1k_1}$. By applying Lemma \ref{Lemma2}, the first term yields a combination of terms of the form $b_j \wedge (1-b_j) \wedge \bigwedge b'_{j,l}$, which is trivial in $K$-theory.

Now let $n>1$,  and assume again w.l.o.g. that all the supports $\A_i$ are contained in $\Nset^n,$ so the $F_i$'s are Taylor polynomials and we can use Lemma \ref{Lemma2}. Fix a solution $\tilde{\alpha}= (\alpha_1 \dots, \alpha_n)$ in $\left(E\setminus\{0\}\right)^n$ for the last $n$ equations. As in the case $n=1$, it is easy to see that we may write equations of the kind
\[x_{l1} = a^{-1} \tilde{F}_l(\tilde{\alpha})\]
where $a$ is equal to a product (possibly empty) of $\alpha_i$.

 We need to consider
\[ z \wedge\bigwedge x_{0j} \wedge \bigwedge x_{1j} \wedge \dots \wedge \bigwedge x_{nj},  \]
which, due to (\ref{pescado}), equals

\[ \sum F_0(\tilde{\alpha}) \wedge \bigwedge_{j=1}^{k_0} x_{0j} \wedge \tilde{F}_1(\tilde{\alpha}) \wedge \bigwedge_{j=2}^{k_1} x_{1j} \wedge \dots  \wedge \tilde{F}_n(\tilde{\alpha}) \wedge \bigwedge_{j=2}^{k_n} x_{nj} \]
\[ + \sum_\eta d_\eta \, {\rm Res}_\eta \wedge
\bigwedge_{j=1}^{k_0} x_{0j} \wedge \tilde{F}_1(\tilde{\alpha}) \wedge \bigwedge_{j=2}^{k_1} x_{1j} \wedge \dots  \wedge \tilde{F}_n(\tilde{\alpha}) \wedge \bigwedge_{j=2}^{k_n} x_{nj} \]
The sum is over all the solutions $\tilde{\alpha}$ of $F_i=0$ for $1 \leq i \leq n$.

For the first term, we apply Lemma \ref{Lemma2} to each of the $n+1$ sets of coefficients and obtain a combination of terms of the form $b_j \wedge (1-b_j) \wedge \bigwedge b'_{j,l}$, which is trivial in $K$-theory.

The terms that correspond to combinations of $\alpha_i  \wedge \bigwedge c_{i,h}$ are zero because we have $n$ different $\alpha_i$'s but $n+1$ terms, which means that some $\alpha_i$ appears twice in the wedge product and then it must be equal to zero.

Consider the second term. For each $\eta$ and $i=1\ldots,n,$ we define $G^\eta_i:= F_i^\eta,$ the restriction of $F_i$ to the facet $\eta$ (see \cite{PS} for a precise definition definition of these polynomials).
As the coefficients of the $G^\eta_i$'s are included in the coefficients of the $F_i$, we can apply induction and obtain the triviality of this term.
This can be done due to the fact that ${\rm Res}_\eta$ is always a resultant associated to a system of dimension $n-1$ or less.
Hence, the symbol is trivial in $K^M_{k+1}(E)\otimes \Qset$ and that proves the claim.
\end{pf*}

\bigskip
\begin{ack}
The authors wish to thank Fernando Rodriguez-Villegas for several suggestions and also  Mart\'{\i}n Sombra for explaining to them his work \cite{som}. They would also like to express their gratitude to the Miller Institute
at UC Berkeley and the Department of Mathematics at the University of Texas at Austin for their support. Matilde Lal\'{\i}n is also grateful to the Harrington fellowship and John Tate for their support. The authors are grateful to the Referee for their suggestions. 
\end{ack}

\end{document}